\documentclass{amsart}

\begin{document}
\parskip10pt
\parindent10pt
\baselineskip15pt

\newtheorem*{theorem}{Theorem}
\newtheorem*{lemma}{Lemma}

\title{A Large Sieve inequality for Euler products}

\author{P.D.T.A. Elliott}
\address{Department of Mathematics, University of Colorado Boulder, Boulder, Colorado 80309-0395 USA}
\email{pdtae@euclid.colorado.edu}

\author{Jonathan Kish}
\address{Department of Mathematics, University of Colorado Boulder, Boulder, Colorado 80309-0395 USA}
\email{jonathan.kish@colorado.edu}

\subjclass[2000]{11N35, 11L20, 11A41}
\keywords{Large Sieve, Euler products, Dirichlet characters}

\begin{abstract}
An inequality of Large Sieve type, efficacious in the analytic treatment of Euler products, is obtained.
\end{abstract}

\maketitle

In this paper we establish an inequality of Large Sieve type that, besides its own interest, lends itself to the study of Dirichlet series with attached Euler products.  Obstacles to sharpenings are discussed.

\begin{theorem}
For each positive real $B$ there is a real $c$ such that
\begin{equation}
\sum\limits_{j=1}^k \max\limits_{y \le x} \max\limits_{\sigma \ge 1, |t| \le D^B} \left|\sum\limits_{D<p \le y} a_p\chi_j(p)p^{-s}\right|^2 \le \left(4L +(k-1)c\right)\sum\limits_{D<p \le x}|a_p|^2p^{-1}, \notag
\end{equation}
with $s = \sigma + it$, $\sigma = \text{Re}(s)$, $L = \sum_{D<p \le x} p^{-1}$, uniformly for $a_p$ in $\mathbb{C}$ and distinct Dirichlet characters $\chi_j \pmod{D}$, $x \ge D \ge 1$.
\end{theorem}
The constant $c$ may be made explicit.

The following result is vital.

\begin{lemma}
Given $B>0$
\begin{equation}
\text{Re} \sum\limits_{w < p \le y} \chi(p)p^{-s} \notag
\end{equation}
is bounded above in terms of $B$ alone, uniformly for $\sigma \ge 1$, $|t| \le D^B$, $y \ge w \ge D$ and all non-principal characters $\chi \pmod{D}$, $D \ge 1$.
\end{lemma}

A proof of this lemma employing analytic properties of Dirichlet $L$-series is given in Elliott \cite{elliottMFoAP6}, Lemmas 1, 4; an elementary proof, via the Selberg sieve and which yields a better dependence upon $B$, is given in Elliott \cite{elliott2002millenium}, Lemma 4; the case $\sigma=1$ by continuity.

\emph{Proof of the Theorem}.  Since the sum $\sum_{D<p \le x} |a_p|p^{-\sigma}$ approaches zero as $\sigma \to \infty$, the innermost maximum may be taken over a bounded rectangle.  In view of the uniformity in $y$, Abel summation allows us to restrict to the case $\sigma=1$.

For reals $t_j$, $y_j$, $|t_j| \le D^B$, $D < y_j \le x$, define
\begin{equation}
\delta_{j,p} = \begin{cases} \chi_j(p)p^{-\frac{1}{2} - it_j} & \text{if} \ D<p \le y_j, \\ 0 & \text{otherwise,} \end{cases} \notag
\end{equation}
$j = 1, \dots, k$, and consider the inequality
\begin{equation}
\sum\limits_{D<p \le x} \left|\sum\limits_{j=1}^k b_j \delta_{j,p} \right|^2 \le \Delta \sum\limits_{j=1}^k |b_j|^2, \notag
\end{equation}
where the $b_j$ are for the moment real and non-negative.  The expanded sum is
\begin{equation}
\sum\limits_{j=1}^k b_j^2L + 2\sum\limits_{1 \le j < \ell \le k} b_jb_\ell \ \text{Re}\sum\limits_{D < p \le x} \chi_j \overline{\chi}_\ell(p)p^{-1-it_j+it_\ell}. \notag
\end{equation}
An appeal to the lemma followed by an application of the Cauchy-Schwarz inequality shows that we may take $\Delta = L+(k-1)c_1$ for a certain $c_1$ depending at most upon $B$.

If now $b_j$ is complex, we represent it as a sum
\begin{equation}
\max(\text{Re} \ b_j,0) + \min(\text{Re} \ b_j,0) + i\max(\text{Im} \ b_j,0)+i\min(\text{Im} \ b_j,0) \notag
\end{equation}
and correspondingly partition the innersum over $j$.  Since the coefficients in each subsum all have the same argument, a second application of the Cauchy-Schwarz inequality allows us to conclude that with $\Delta = 4\left(L+c_1(k-1)\right)$ the above inequality holds for all complex $b_j$.

Dualising:
\begin{equation}
\sum\limits_{j=1}^k \left|\sum\limits_{D<p \le x} a_p \delta_{j,p} \right|^2 \le 4(L+c_1k)\sum\limits_{d<p \le x} |a_p|^2 \notag
\end{equation}
for all complex $a_p$, c.f. Elliott \cite{elliott1971inequalitieslargesieve}.

Replacing $a_p$ by $a_pp^{-\frac{1}{2}}$ completes the proof.

\emph{Remarks}.  Presumably the inequality in the Theorem is valid with the coefficient $4$ replaced by~$1$.  With the present argument that appears to require the sum in the Lemma, with $\sigma=1$, $w=D$, to be uniformly bounded not only above but also below.  Such a bound seems currently out of reach.  It would, in particular, guarantee a lower bound $L(1,\chi) \ge c_2(\log D)^{-1}$ for quadratic characters $(\text{mod} \ D)$ and eliminate Siegel zeros.

Without an adjustment to the term $(k-1)c$ the restriction $D<p$ in the sums over the primes cannot be altogether removed:

Application of the Large Sieve, c.f. Elliott \cite{elliott1969onthesizeofL1chi}, \cite{elliott1970distributionquadraticclassnumber}, \cite{Elliott1980}, shows that if $(\log x)^{20} \le y \le x^2$, then with at most $O(x^{\frac{7}{8}})$ exceptions
\begin{equation}
\sum\limits_{p > y} \chi(p)p^{-1} \ll y^{-\frac{1}{10}} \notag
\end{equation}
is satisfied by the primitive characters to moduli $D$ not exceeding $x$.  With a slight increase in the number of exceptions, the method of \cite{elliott1969onthesizeofL1chi} guarantees a similar result that is uniform in $y$.

We may identify Dirichlet characters $\chi$ of order $m$ to prime moduli $q$ with $m^{th}$-power residue symbols and view them in terms of characters on ideal class groups, as in Elliott \cite{elliott1970meanvalueoffofp}, where appropriate references to works of Eisenstein, Landau, Furtw\"angler, Artin and Hasse may be found.

Employing the uniform distribution of prime ideals in ideal classes, in particular Fogels' generalisation of Linnik's theorem on the size of the least prime in a rational residue class, c.f. Elliott~ \cite{elliott1967problemoferdos}, \cite{elliott1968notesonkthpowerresidues}, \cite{elliott1970distributionofpowerresidues}, Fogels \cite{fogels1961distributionofprimeideals}, one may arrange an infinitude of moduli $q$ for which $\chi(p)$, with $(p,m)=1$, $p$ up to a certain constant multiple of $\log q$, may be given individually any value available to the character.

As an example, if
\begin{equation}
\min\limits_{1 \le r \le m} \cos\left(2\pi rm^{-1}\right) \le \beta \le 1, \notag
\end{equation}
then by choosing the successive values of $\chi(p)$ to be complex conjugates we may arrange that
\begin{equation}
\sum\limits_{p \le q} \chi(p)p^{-1} = \beta \log\log\log q + O(1). \notag
\end{equation}
Moreover, with $\beta = \pm 1$, separately, the estimate may be required to hold for every character of order $m$ or even order $m$, respectively.

Via the construction of finite probability spaces, these methods allow the successful study of the values of series $L(s,\chi)$ formed with Dirichlet characters of order $m$ to prime moduli provided $\sigma > 1-c_3>\frac{1}{2}$, reaching part-way into the critical strip, c.f. Elliott \cite{elliott1973distributionquadraticLseries}.

More generally, taking imprimitive characters into account, for almost all moduli $D$, in a strong quantitative sense, the sums in the Lemma are indeed bounded below and the inequality of the Theorem valid with $4$ replaced by $1$.

Variant inequalities also allow the constant $4$ to be reduced.  For example, replacing $\chi_j(p)p^{-s}$ by $\text{Re}\left( \chi_j(p)p^{-s}\right)$ we may replace $4\left(L+(k-1)c\right)$ by $2\left(L+kc\right)$.  Note that summands corresponding to a complex character $\chi_j$ may then appear twice in the bounded sum.

\bibliographystyle{amsplain}
\bibliography{c:/Me/Boulder/Mathbib}

\end{document}